\input amstex
\documentstyle{amsppt}
%\magnification 1200
\NoBlackBoxes
\input epsf

                     %     MACROS
  % to use in \if{}\else{}\fi
  %

%

\def\Z {\bold Z }

\def\C {\bold C }

\def\R {\bold R }

\def\<{\langle}
\def\>{\rangle}
\def\U{\sqcup }

\def\Card{\operatorname{Card}}
\def\Binom{\binom}
\def\pr{\operatorname{pr}}

\def\wsim{\asymp}
\def\esim{\underset{e}\to\sim}
\def\wesim{\underset{e}\to\asymp}
\def\ewsim{\wesim}

                          % References

\def\refAndrews      {1} 
\def\refArnold       {2} %Arnold: Integral polygons
\def\refDK	     {3} %Degtyarev Kharlamov (UMN)
\def\refGH           {4} %Griffits, Harris
\def\refGKZ          {5} %Gelfand-Kapranov-Zelevinski
\def\refHaas         {6} %Haas
\def\refHarnack      {7} %Harnack
\def\refHP           {8} %Harari-Palmer
\def\refItenberg     {9} %Itenberg
\def\refKS           {10} %Konyagin-Sevastianov
\def\refOrevkovJCTA  {11} %Orevkov:JCTA
\def\refOrevkovVol   {12} %Orevkov:Volume of discriminant
\def\refOtter        {13} %Otter
\def\refPolya	     {14} %Polya
\def\refThom         {15}
\def\refVershik      {16}
\def\refViro         {17} %Viro: Glueing
\def\refTree         {18} %Math.Encycl

                           % FIGURES

                           %

                      %     DISPLAYED FORMULAS

                      %

\def\mRec          {1}
\def\NRec          {2}
\def\eqOtter       {3}
\def\eqConst       {4}

\rightheadtext{ %Growth rate
Asymptotic growth
of the number of classes }
%of real curves }
\topmatter

%-----

\title          %Growth rate
                Asymptotic growth of the number of classes of real
                plane algebraic curves when the degree increases
\endtitle

%-----

\author          S.Yu.Orevkov, V.M.Kharlamov
\endauthor

%-----

\address         Steklov Math.Inst. (Moscow, Russia) \endaddress

%-----

\address         Universit\'e Paul Sabatier (Toulouse, France)
\endaddress

%-----

\address         
IRMA (CNRS) et Universit\'e Louis Pasteur, Strasbourg, France
\endaddress

%-----

\endtopmatter

%\hskip1.15in EXACT ANSWERS are great when we can find them.
%
%\hskip1in R.L.Graham, D.E.Knuth, O.Patashnik, {\it Concrete
%Mathematics}

%============================================= epigraph
%
\font\eightrus=wncyr8
\font\eightirus=wncyi8

\hskip2in {\eightrus Chto my znaem o lise? Nichego! I to --- ne vse.}

\hskip3in {\eightrus B.~Zahoder.} {\eightirus Mohnataya azbuka}
%================================================================

\document

\head Introduction
\endhead

Plane algebraic curves and, more generally, projective
algebraic hypersurfaces is the subject of the 16th Hilbert problem.
They are naturally organized in families which are the spaces
of homogeneous polynomials. These spaces are numbered by the number
$n+1\ge 2$ of homogeneous
variables and the degree $d\ge 1$
of polynomials. Each one is a projective space
$P^N$ with $N=N(n,d)-1$, $N(n,d)=\binom {n+d} d$
being the number of
coefficients in a generic homogenous polynomial of degree $d$
in $n+1$ variables. One associates with it
the universal hypersurface
$\Gamma\subset P^N\times P^n=\{(p,x)\,\vert\, p(x)=0\}$,
which is a nonsingular variety, and the projection
$\pr :\Gamma\to P^N$. The critical locus $D\subset P^N$
of $\pr$
is the so-called {\it discriminant} hypersurface:
$p\in D$ if and only if the hypersurface
$p=0$ in $P^n$ is singular, and the
induced over $P^N\setminus D$ family
$\Gamma^0=\pr^{-1}(P^N\setminus D)\to
P^N\setminus D$ is a deformation family:
$\pr\vert_{\Gamma^0}$ is a proper submersion.

Over $\C$,
the space of complex points of $P^N\setminus D$ is
connected, which implies that
all the nonsingular
hypersurfaces of same degree (and same dimension) are diffeomorphic.
Over $\R$,
the situation is completely different: as soon as $d\ge 2$,
$P^N_\R\setminus D_\R$ (here and further we denote the real
point set of a real variety $X$ by $X_\R$)
is disconnected.
There is a tradition to call
{\it rigid isotopy} or {\it real deformation}
a continuous path in $P^N_\R\setminus D_\R$;
respectively, two real
hypersurfaces are {\it rigid isotopic} or {\it 
real deformation   equivalent} if they are connected by a rigid isotopy. 
The
simplest
invariant of the rigid isotopy is the topology of the real
part of the hypersurface. 
%It is sufficient 
%to prove 
By means of it one proves easily
the above disconnectedness in any dimension and
any degree $\ge 2$.

It is worth to be noted that over $\R$,
the definition of $D$
can be interpreted in two different ways: $p\in D$ if
$p=0$ has a real singular point, or $p\in D$ if $p=0$
has a complex singular point. We prefer the latter,
commonly used, definition,
which is equivalent to considering
the real discriminant as the real point set
$D_\R\subset P^N_\R$ of the complex one
(and which gives a discriminant defined over $\Z$). In fact,
for our purposes, see below,
%counting the number of the connected components
%of the complement
there is no difference: the
complement of the other one is
included into $P^N_\R\setminus D_\R$,
and the inclusion establishes a one to one correspondence
between the connected components, since the difference
between the two discriminants is of dimension $\le N-2$.

In his 16th problem Hilbert asks about "the number, 
form, and position of the sheets" of a nonsingular real
algebraic hypersurface of given dimension and
degree (in fact, he is mentioning
only the plane curves and the surfaces in 3-space).
In its extended form the problem can be understood as the
study of the deformation equivalence and the deformation invariants.

  Let $I(n,d)$ be the number of isotopy types of pairs
  $(P^n_\R, X_\R)$ where $X_\R$ is a nonsingular hypersurface
  in $P^n_\R$ of degree $d$ and let $D(n,d)$ be the number
  of deformation classes (rigid isotopy classes).
  It is clear that $I(n,d)\le D(n,d)<\infty$.
  In what follows we often abbreviate
  $I(2,d)$ and $D(2,d)$ to $I_d$ and $D_d$. 
(If $n=2$
the isotopy equivalence coincides with equivalence by
homemorphisms of $P^2_\R$ and the isotopy classes have 
a simple encoding, see below.)

The numbers (and the corresponding classes) $I(n,d)$ 
and $D(n,d)$ are known for some values of $(n,d)$:
%-------------------------------------------------------------
$$
        I(1,d)=D(1,d)=[\frac{d}2]+1,
        \quad
        I(n,1)=D(n,1)=1,
        \quad   
        I(n,2)=D(n,2)=n,
$$

and
\footnote{
	The bound $I(3,5)\ge 3000$ was recently announces 
	by B.~Chevallier (private communication).
} 

$$
\xalignat{9}
             d:&& 1 && 2 && 3 &&  4 &&   5 &&  6 &&   7 &&    8&     \\
              :&                                                     \\
        I(2,d):&& 1 && 2 && 2 &&  6 &&   8 && 56 && 121 && \ge 2500& \\
        I(3,d):&& 1 && 3 && 5 &&111 &&\ge 3000 &                     \\
              :&                                                     \\
 D(2,d)-I(2,d):&& 0 && 0 && 0 &&  0 &&   1 &&  8 &&\ge 570 &         \\
 D(3,d)-I(3,d):&& 0 && 0 && 0 && 55 &&   &                           \\
\endxalignat
$$

For higher degrees and dimensions,
very few is known. To our knowledge, 
even the asymptotic behavior
of $I(n,d)$ and $D(n,d)$
were never studied
systematically, and in this paper we are making
an attempt to formulate the principal questions
and to give some answers concerning
the asymptotics up to weak exponential equivalence
(see the definition below). They are more advanced in
the case of plane curves, 
where, in particular, we give the asymptotics
for $I(2,d)$ and for the number of 
isotopy classes of maximal curves
realizable by $T$-curves
and show that none of the known restrictions 
is asymptotically valuable.
%and show that the restrictions coming
%from the B\'ezout theorem
%are the most valuable.

{\bf Notation.}
Beside the usual
notion of the {\it equivalency}
 of sequences
($a_m\sim b_m$ if $\lim_{m\to\infty} a_m/b_m\,=1$),
we shall often use 
{\it weak equivalency relation}
%)
$$
        a_m \wsim b_m \qquad\text{ if }\quad
	a_m=O(b_m)\,\, \text{ and } \,\, b_m=O(a_m),
$$
{\it exponential equivalency}
($a_m\esim b_m$ if $\log a_m \sim\log b_m$) and
{\it exponential weak equivalency}
($a_m\ewsim b_m$ if $\log a_m \wsim\log b_m$).

\example{ Example }
$m! \ewsim m^m$ (by Stirling formula); $(m^2)^{m^2} \ewsim m^{m^2}$.
\endexample

The real point set $A_\R$ of a real nonsingular plane curve
$A\subset P^2$ is
%either
a compact nonsingular $1$-dimensional
smooth submanifold of $P^2_\R$, or empty.
Any compact nonsingular $1$-dimensional submanifold $C$
of
$P^2_\R$ consists of
a finite number of disjoint embedded circles.
We call $\Z/2$-degree of $C$ the element
realized by $C$ in
$H_1(P^2_\R;\Z/2)=\Z/2$.
%If the $C$ realizes zero in $H_1(P^2_\R;\Z/2)$,
If it is zero,
all the circles are two-sided. Otherwise, the curve has one and only
one one-sided component.  The two-sided circles are called {\it
ovals}. In the algebraic case, $C=A_\R$, the curve realizes zero
if and only if its degree is even.

Any oval decomposes $P^2_\R$ in two parts: the internal one,
homeomorphic to a disc, and the external one, homeomorphic to
a cross-cup (Moebius band).
Thus, there appears a natural partial order (tree structure):
$a>b$ if the oval $a$ contains $b$ in its interior.
This order is invariant under isotopy
(which, in our case, is equivalent to invariance under
homeomorphisms $P^2_\R\to P^2_\R$),
and, conversely, this order and $\Z/2$-degree determine the
curve up to isotopy.

Following tradition, a set of disjoint embedded circles
considered up to an isotopy in $P^2_\R$
will be called a {\it real scheme} (of an oval arrangement).
Any real scheme can be realized by an algebraic curve.
A real scheme is said of degree $d$ if it
can be realized by nonsingular algebraic curves of
degree $d$.

We represent a real scheme by the rooted tree with $l+1$
vertices, where $l$ is the number of ovals (the root is an
additional, the greatest, element; we introduce it to be in
accordance with the terminology established in the combinatorics of
trees), and, as well, by nested
parentheses. In the latter notation,
the now traditional 
rule is such that
$\<1\<A\>\>$ denotes a scheme $\<A\>$ with one additional oval
containing $\<A\>$ inside it and $\<A\U B\>$ the disjoint sum of the
schemes $\<A\>$ and $\<B\>$. (We
omit $\Z/2$-degree, since in what follows it never leads to a
confusion.)

Note that in the algebraic case,
according to the Harnack inequality,
the number of ovals
is bounded by $\frac12(d-1)(d-2)+1$ if $d$ is even,
and by $\frac12(d-1)(d-2)$ if $d$ is odd. The
real schemes of degree $d$ with such a maximal number
of ovals are called $M$-schemes and the
corresponding curves are called $M$-curves.  They exist for any
degree.

{\bf Acknowledgements.}
We would like to thank
I.~Itenberg, S.~Fomin, M.~Lifschitz, E.~Shustin, A.~Vershik,
with whom we discussed these and related problems.
This work was partially supported by Program "Arc en Ciel 2000".

\head        1. Real plane curves
\endhead

\subhead     1.1. Curves up to isotopy
\endsubhead
The following statement describes the rough
asymptotics of the number of curves considered up to
isotopy.

\proclaim{ Proposition 1 }    $I_d \ewsim \exp(d^2).$
\endproclaim

\demo{ Proof }
{\it Upper bound. }
By the Harnack inequality the number of ovals is bounded
by $\frac12(d-1)(d-2)+1$. Hence, $I_d$
does not exceed the number of rooted trees with
$\frac12(d-1)(d-2)+2$ vertices. It is known (see, e.g.,
\cite{\refTree}), that the number of rooted trees with $n$ vertices
is bounded by $C^n$ for some constant $C>0$. (For the precise value
of the constant see \cite{\refOtter}; the first digits are pointed in
Section 1.4 below).

{\it Lower bound.}
Applying Viro's gluing method, we inductively construct
for each $k=0,1,2,\dots$ a finite family $\Cal C_k$
of pairwise non-isotopic nonsingular
real plane curves of degree  $d_k$, where
$$
         d_0=6,\qquad d_{k+1} =2d_k+6
                                                \eqno(\mRec)
$$
in such a way that the real point set
$c_\R$ of each curve $c$ from $\Cal C_k$
is contained in the affine plane $\R^2$.

To %define
   %build
   construct
$\Cal C_0$, we pick
a representative from each isotopy
class of curves of degree $d_0=6$
(they all can be chosen in $\R^2$ because, as is known,
each one is realizable by a perturbation of three
suitable ellipses).

As soon as the family $\Cal C_k$ is constructed,
associate with every curve $c\in C_k$
a curve which is the union of a round
circle lying in the positive quadrant and the image $f(c)$ of $c$
under an affine linear transformation of the plane such that $f(c_\R)$
is mapped inside the circle.
Denote by $\Cal C'_k$
the set of the
obtained curves of degree $d_k+2$. Then,
for any non-ordered $4$-tuple $c'_1,\dots,c'_4\in \Cal C'_k$,
construct a curve by Viro's method
dividing the triangle
$\Delta=[(0,0),(0,2d_k+6),(2d_k+6,0)]$ 
in four triangles
with the side $d_k+2$  separated one from another by two ribbons of
width $1$ and one ribbon of width $2$
(three of these four triangles have a vertex
and two sides on $\partial\Delta$;
the ribbons are necessary to glue the neighboring charts
between them). Let $N_k=\Card \Cal C_k$.  
It is clear from the above construction that
$$
   N_{k+1}\ge N_k^4/24          \eqno(\NRec)
$$
We get by induction from (\mRec) and (\NRec) that
$$
        d_k =2^k\cdot(d_0+6)-6,
        \qquad \log N_k\ge 4^k\cdot\big(\log
     N_0-{\log24\over3}\big)+ {\log24\over3}.\qed
$$

\enddemo

%%%%%%%%%%%%%%%%%%%%%%%%%%%%%%%%%%%%%%%%%%%%%%%%%%%%%%%%%%%%%%
%%%%%%%%%%%%%%%%%%%%%%%%%%%%%%%%%%%%%%%%%%%%%%%%%%%%%%%%%%%%%%
%%%%%%%%%%%%%%%%%%%%%%%%%%%%%%%%%%%%%%%%%%%%%%%%%%%%%%%%%%%%%%
%%%%%%%%%%%%%%%%%%%%%%%%%%%%%%%%%%%%%%%%%%%%%%%%%%%%%%%%%%%%%%

\subhead     1.2. Curves up to deformation
\endsubhead
The following proposition gives a first information
on the growth of the number of curves considered up to
rigid isotopy.

\proclaim {Proposition 2}
There are constants $c_1, c_2$ such that  for any $d$
$$
   c_1 d^2\le \log D_d\le c_2 d^2\log d.
$$
\endproclaim

\demo{Proof}
The lower bound follows from Proposition 1.
By Poincar\'e-Lefschetz duality, to prove the upper bound
it is sufficient to get the same bound for the total
Betti number of the real point set $D_\R$
of the discriminant hypersurface of singular curves of degree $d$.
Such a bound is given, e.g., by the Smith-Thom inequality:
indeed, the degree $h$ of the discriminant hypersurface
is $3(d-1)^2$ and the leading term of this bound
(polynomial in $h$ and exponential in $N$)
is $h^N$ where
$N=N(2,d)-1= {\binom {d+2} d} -1$
(see \cite{\refThom}).
\qed
\enddemo

By a $T$-curve of degree $d$
we mean a (real nonsingular) plane curve
obtained by Viro's gluing method applied to
a convex triangulation of the Newton triangle
$[(0,0), (0,d), (d,0)]$
in primitive (i.e., area $\frac12$)
triangles with integral vertices
(see \cite{\refViro}, \cite{\refItenberg}, \cite{\refHaas}).
Denote by $D^T_d$ the number of deformation classes
realized by $T$-curves of degree $d$.

\proclaim{Proposition 3}
There is a constant $c>0$ such that
$\log D^T_d\le cd^2$ for any $d$.
\endproclaim

\demo{Proof}
As it follows from Gelfand-Kapranov-Zelvinskii description
of the secondary polytopes (see \cite{\refGKZ}),
a convex triangulation of
a (convex) polygon in primitive integral triangles is determined
by the multiplicities of the vertices of the triangulation.
The total number of vertices is $\frac12(d+1)(d+2)$
and the number of edges is $\frac32d(d+1)$. Thus,
the number of convex triangulations is bounded from above
by the number of decompositions of
$3d(d+1)-(d+1)(d+2)+3$ in $\frac12(d+1)(d+2)$
integral positive
summands, which is the binomial coefficient
$$
  \binom {2d^2} {\frac12(d+1)(d+2)-1} \;\esim\; c^{\frac12d^2}
  \qquad\text{with}
  \quad c=\frac{4^4}{3^3}.
$$
The choice of signs in the vertices
can only change the final constant $c$.
\qed

\enddemo

%%%%%%%%%%%%%%%%%%%%%%%%%%%%%%%%%%%%%%%%%%%%%%%%%%%%%%%%%%%%%%
%%%%%%%%%%%%%%%%%%%%%%%%%%%%%%%%%%%%%%%%%%%%%%%%%%%%%%%%%%%%%%
%%%%%%%%%%%%%%%%%%%%%%%%%%%%%%%%%%%%%%%%%%%%%%%%%%%%%%%%%%%%%%
%%%%%%%%%%%%%%%%%%%%%%%%%%%%%%%%%%%%%%%%%%%%%%%%%%%%%%%%%%%%%%

\subhead    1.3. $M$-curves
\endsubhead
The aim of this section is to recover the asymptotic growth
of the number of $M$-schemes which are realizable by $T$-curves.

It is worth noting that similar to $I_d$ the
whole number of pairwise non-isotopic $T$-curves of degree $d$
is exponential weak equivalent
to $\exp(d^2)$ (one can slightly modify the construction in the proof
of the lower bound of Proposition 1).
Combined with the upper bound from Proposition 3,
this shows that the number of rigid
isotopy types of degree $d$ realizable by $T$-curves
is also exponential weak equivalent to $\exp(d^2)$.

\proclaim{ Proposition 4} {\rm(see [\refItenberg])}
If an $M$-scheme of degree 
$d$ is realizable by a $T$-curve
then it satisfies the following condition:
\roster
\item "(*)"
            The number of ovals $O$ such that
            there exist ovals $O'$ and $O''$ such that
            $O''<O'<O$, does not exceed 
$\frac32 m$.\qed
\endroster
\endproclaim

\proclaim{ Proposition 5}
Let $T^*_d$ be the number of different schemes which
satisfy the condition {\rm(*)} and have at most 
$\frac12(d-1)(d-2)+1$ ovals.
Then  $T^*_d\wesim\exp(d^{3/2})$. 
\endproclaim

\demo{ Proof }
Let $P(m)$ 
denotes the set of all non-ordered partitions of
$m$, i.e. 
$P(m)=\{(\lambda_1,\dots,\lambda_k)\,|\,\lambda_1\ge\dots\ge \lambda_k>0,\;
\lambda_1+\dots+\lambda_k=m\}$, and let 
$p(m)=\Card P(m)$.
It is known that $p(m)\wesim \exp(\sqrt m)$ (see \cite{\refAndrews}).
For $\lambda =(\lambda_1,\dots,\lambda_k)\in P(m)$, let
$S(\lambda)$ denotes
the real scheme of oval arrangement
$1\<\lambda_1-1\>\U\dots\U 1\<\lambda_k-1\>$.

{\it Lower bound.}
For any ordered collection
$(\lambda^{(1)},\dots,\lambda^{(k)})\in P(d)^k$
where $k=[3d/2]$, the scheme
$ S_1 \U 1\< S_2 \U 1\<\dots\U 1\< S_k\>\dots\>$
 where
$S_j=S(\lambda^{(j)})$,
satisfies the condition (*).
Hence,
$T^*_d\ge p(d)^k\wesim \exp(d\sqrt d)$.

{\it Upper bound.}
Let $S$ be a real $M$-scheme of degree 
$d$ satisfying the condition (*).
Denote by $S^*$ the subscheme of $S$ consisting of
the ovals which are the outermost
ovals of nests of the depth
$\ge 3$ (by the condition
(*), we have $\Card S^* \le \frac{3d}2$).
Let $T(S^*)$ be the rooted tree of $S^*$.

Let us choose a representative in each class of isomorphic trees
and let us fix (arbitrarily) a numbering of its vertices.
Let $v_1,\dots,v_k$ ($k\le \frac{3d}2$) be the ovals of $S^*$ numbered
in accordance with the chosen order of the vertices of $T(S^*)$.
For any $j$ with $1\le j\le k$, let us consider all the ovals of
$S\setminus S^*$ lying inside $v_j$ which are not separated from
$v_j$ by other ovals of $S^*$.
These ovals constitute a real subscheme of the form 
$S(\lambda^{(j)})$ for some partition
$\lambda^{(j)} \in P(m_j)$, and $m_1+\dots+m_k=l(d)-k$ where
$l(d)\sim d^2/2$ is the number of ovals in an $M$-scheme.
Thus,
$$
        T^*_d\le \sum_{k=1}^{[3d/2]} T_k
        \sum_{m_1+\dots+m_k=l(d)-k} p(m_1)\dots p(m_k)
$$
where $T_k$ is the number of all rooted trees with $k$ vertices.
For some constant $a$ we have $\log p(m_j)<a\sqrt{m_j}$.
Hence, by the Cauchy-Buniakowski inequality, one can uniformly
estimate each summand in the interior sum:
$$
  \log\big(p(m_1)\dots p(m_k)\big)<a\sum\sqrt{m_j}\le
  a\sqrt{m_1+\dots+m_k}\sqrt{k}<a\sqrt{l(d)\cdot3d/2},
$$
the number of the summands being
$\Binom{l(d)}{k}< \Binom{d^2}{3d/2}\ewsim d^d$
which implies
$$
  T^*_d<(3d/2)\,T_d \Binom{d^2}{3d/2} \exp(b d^{3/2})\ewsim \exp(d^{3/2}).
\qed
$$
\enddemo

Let $I^{TM}_d$
denotes the number of $M$-schemes of degree $d$
which are realizable by $T$-curves.

\proclaim{ Proposition 6 }
$I^{TM}_d\wesim \exp( d^{3/2})$.
\endproclaim

\demo{ Proof }
{\it Upper bound.} Follows from Propositions 4 and 5.

{\it Lower bound.}
Let $d$ be a positive integer.
   %Denote
   %Pose
   Set
$k=[\frac16 d]-1$ and
denote the triangle 
$[(0,0),(d,0),(0,d)]$
by $\Delta_d\subset \R^2$.
For any collection of partitions
$(\lambda_0,\dots,\lambda_{k-1})$,
$\lambda_j\in P(2(k-j))$ we shall construct an 
$M$-curve by Haas' method \cite{\refHaas}.

Let us consider the following points in $\Delta_d$:
$$
\split
       &A_j=(5j,0),\;\; j=0,\dots,k+1;\quad
        B_0=(0,d),\;\; B_j=(5j+1,6(k-j)-2),\;\; j\ge1;
                                                       \\
       &C_j=(6k-j,0),\;\;
        D_j=(5j+1,0),\;\;
        E_j=(5j+2,6(k-j)-4),\;\; j=0,\dots,k.
\endsplit
$$
Let us cut 
$\Delta_d$
by distinct two-segment broken lines
$$
        A_jB_jC_j\;\; (j=1,\dots,k)
        \qquad\text{and}\qquad
        D_jE_jA_{j+1}\;\; (j=0,\dots,k).
$$
Now let us cut each triangle $D_jE_jA_{j+1}$
($j=0,\dots,k-1$) as follows.
Let $\lambda_j=(\lambda_j^{(1)},\lambda_j^{(2)},\dots\,)$.
Set $s_j^{(\nu)}=\lambda_j^{(1)}+\dots+\lambda_j^{(\nu)}$,
$y_j^{(\nu)}=2+2\cdot[3s_j^{(\nu)}/2]$, and 
$E_j^{(\nu)}=(5j+2,y_j^{(\nu)})$.
Let us cut each triangle $D_jE_jA_{j+1}$
into domains $Z_j^{(1)},Z_j^{(2)},\dots$
by the broken lines $D_jE_j^{(1)}A_{j+1}$,
$D_jE_j^{(2)}A_{j+1},\dots\;$.
Then the interior of the segment which is cut by the domain
$Z_j^{(\nu)}$ from the line $x=5j+3$, contains
$\lambda_j^{(\nu)}$
integral points with an even $y$-coordinate
(totally, the segment which is cut by the triangle $D_jE_jA_{j+1}$ from
the line $x=5j+3$, contains exactly $2(k-j)$ points with an even
$y$-coordinate).

Thus, we have cut $\Delta_d$ into $k+s+1$ polygons where
$s$ is the total number of elements of the partitions 
$\lambda_1,\lambda_2,\dots$.
It is not difficult to check that any side of any of the polygons
contains no other integral points except the ends, and one of the ends
is on the ends lies on a side of $\Delta_d$.
Following \cite{\refHaas}, we shall call these polygons {\it zones}.
Let us attach a sign to each zone so that the signs of adjacent zones
are opposite.
Then, extend the zone decomposition to a convex
triangulation and define a sign distribution
$\delta:\Delta_d\cap\Z^2\to\{\pm1\}$ setting
$\delta(x,y) =(-1)^{xy}\zeta^y$ where $\zeta$ is the sign of the zone
which contains the point $(x,y)$
(since all integral points on the zone boundaries have even
$y$-coordinates,
the definition is coherent on the intersections of the zones).

By Haas' theorem \cite{\refHaas},
the corresponding $T$-curve is an $M$-curve and it is not difficult
to check that its real scheme has form
$$
        a\U 1\<a_0\U S_0 \U1\<a_1\U S_1 \U1\<  \;\;\;\;\dots\;\;
        \U1\<a_{k-1}\U S_{k-1} \U1\<1\>\>\dots\>\>\>
$$
for some integers $a,a_0,a_1,\dots\,$,
where $S_j=1\<\lambda^{(1)}_j\>\U1\<\lambda^{(2)}_j\>\U\dots$ for
$\lambda_j=(\lambda^{(1)}_j,\lambda^{(2)}_j,\dots\,)$.
The nest of the depth $k+1$ corresponds to the broken lines
$A_jB_jC_j$, and the subschemes $S_j$ correspond to cuttings of
the triangles $D_jE_jA_{j+1}$. The nests
$1\<\lambda^{(1)}_j\>$, $1\<\lambda^{(3)}_j\>,\dots\,$ 
are above the axis $y=0$,
and the nests $1\<\lambda^{(2)}_j\>$, $1\<\lambda^{(4)}_j\>,\dots\,$ 
are beneath the axis $y=0$.

Thus, for any $d$, we constructed
$$
        \prod_{j=1}^k p(2j) \wesim \exp(\sum_{j=1}^k \sqrt{j})
        \wesim 
        \exp(d^{\frac32})
$$
different $M$-schemes. \qed
\enddemo

%%%%%%%%%%%%%%%%%%%%%%%%%%%%%%%%%%%%%%%%%%%%%%%%%%%%%%%%%%%%%%
%%%%%%%%%%%%%%%%%%%%%%%%%%%%%%%%%%%%%%%%%%%%%%%%%%%%%%%%%%%%%%
%%%%%%%%%%%%%%%%%%%%%%%%%%%%%%%%%%%%%%%%%%%%%%%%%%%%%%%%%%%%%%
%%%%%%%%%%%%%%%%%%%%%%%%%%%%%%%%%%%%%%%%%%%%%%%%%%%%%%%%%%%%%%

\subhead     1.4. A coefficient of 
                  $d^2$ in the asymptotic of the
                  number of real schemes
\endsubhead
In Proposition 1, we used the number of rooted trees as an upper
bound for $I_d$. Due to Otter \cite{\refOtter} 
(see also \cite{\refHP; Section 9.5}),
the following exponential equivalence for the number 
$T_n$ of rooted trees holds
$$
        T_n\esim C^n,
         \qquad C=2.95576                           \eqno(\eqOtter)
$$
Hence,
$$
        \log I_d  \le  A d^2 + o(d^2),
        \qquad  A=(\log C)/2.
                                                 \eqno(\eqConst)
$$
It happens, that none of the known restrictions for the arrangement
of ovals allows us to reduce the coefficient $A$ in the estimate
(\eqConst).
The results of this sections do not serve to improve the upper
bound for $I_d$ (which is one of the goals of
the paper) but in contrary, they show that something is
useless for this purpose. 
This is why we do not give proofs in all detail.

The only known restrictions which could be asymptotically valuable,
are the restrictions coming from B\'ezout's theorem 
(for instance, the absence of nests deeper than $d/2$).
Other restrictions (Petrovsky inequality, Arnold inequality etc.)
provide corrections in (\eqConst) of order $O(1)$ or $O(\log d)$).

To give a strict sense to the statement of non-improvability of
the estimate (\eqConst), let us formalize the notion of a restriction
coming from B\'ezout's theorem.
Let us say that {\it a real scheme $S$ satisfies B\'ezout's theorem
for the degree $d$} if for any triple of positive integers
$(k,d',c)$, $d'<d$ satisfying the condition
\roster
  \item"(**)" {\it Through any
  $k$ points in $P^2_\R$ in general position, there exists a real curve
  of the degree $d'$ which has at most $c$ connected components}
\endroster
and for any choice of $k$ points there exist $c$ smoothly embedded
circles passing through the chosen points which intersect $S$ at
most in $dd'$ points.

\remark{ Remark }
As it follows from Harnack's proof of Harnack's inequality
\cite{\refHarnack}, any real
scheme satisfying B\'ezout's theorem for a degree $d$ has at most
$\frac12(d-1)(d-2)+1$ connected components.
\endremark

Let 
$B_d$ denotes the number of different schemes satisfying B\'ezout's
theorem for the degree $d$.

\proclaim{ Proposition 7}
$\log B_d \sim A d^2$ where $A$ is the
constant in (\eqOtter) and (\eqConst).
\endproclaim

A proof follows immediately from Lemmas 8 and 9 below.

\definition{ Definition } A {\it nest} is a sequence of
disjoint ovals $v_1,\dots,v_p$ in $P^2_\R$ such that
$v_1>\dots>v_p$. The number of ovals $p$ is
called the {\it depth} of the nest.
The maximal depth of a nest of a real scheme (of a rooted tree)
$S$ is called the {\it depth} of $S$.
\enddefinition

A triple $(k,d',c)$ with $d'<d$ satisfies the condition (**)
only if (conjecturally, if and only if)
               $k\le\frac12 d'(d'+3)$ and $k\le3d'+c-2$,
cf. \cite{\refDK}.
Using this fact, one can prove the following statement.

\proclaim{ Lemma 8 }
        Let $\alpha$ be a fixed number such that $0<\alpha<\frac16$.
        Let $d$ be a sufficiently large integer.
        Suppose that a real scheme $S$ has
        at most $l(d)=\frac12\big((d-1)(d-2)+(1+(-1)^d)\big)$
        ovals and $S$ is of depth $<[\alpha d]$.
        Then $S$ satisfies B\'ezout's theorem for degree $d$.\qed
\endproclaim

Let us denote by 
$T_n^{[h]}$ the number of rooted trees 
of depth $\le h$ with $n$ vertices.

\proclaim{ Lemma 9 }
Let $h_n$ be a sequence such that 
$$
    \lim_{n\to\infty}{\log n\over h_n}=0. \qquad
        \text{Then } \;T_n^{[h_n]} \esim T_n \esim C^n.
$$
\endproclaim

\demo{ Proof } Let us number arbitrarily the vertices of each tree.
Pick a rooted tree $T$ of depth $> h$
with $n$ vertices. Choose the highest vertex $v$.
If it is not unique, choose the one whose number is minimal.
Consider the branch of the depth $h-1$ which contains $v$
(by definition, a {\it branch} of a tree is one of two components
obtained by deleting an edge). Cutting this branch, we obtain
another rooted tree. If its depth is still $>h$, 
we iterate this process. Since each time we cut out
$\ge h$ vertices, we perform not
more than $n/h$ iteration. As the result, 
we obtain $\le n/h$ rooted trees of depth $\le h-1$.
Connecting all of them to a new common root, we obtain a
rooted tree with $n+1$ verices of depth $\le h$.
To reconstruct $T$, it suffices to indicate the ends of edges
which were deleted at each step. Thus,  
$T_{n} \le T_{n+1}^{[h]}\cdot n^\frac{n}{h}$.
%\mnote{Kh: $n/h$ instead of $2h$, is not it? 
%       Or: $2n/h$ because each edge has 2 ends}          
\qed
\enddemo

\remark{ Remarks }
{\bf 1.} Our proof of Proposition 1 provides the lower bounds
$$
\split
        &I_d \ge C_6^{d^2},\quad C_6=1.02081\qquad(d_0=6,\,N_0=56)    \\
        &I_d \ge C_8^{d^2},\quad C_8=1.03511\qquad(d_0=8,\,N_0\ge2500).
\endsplit
$$
We see that both constants $C_6$ and $C_8$ are very far from $C$.
\smallskip

{\bf 2.} The numbers $T_n^{[h]}$ can be computed by recurrent formula
for the generating functions \cite{\refPolya}
(see also \cite{\refHP; Section 3.1}):
$$
    T^{[h+1]}(x)=x\exp\Big(\sum_{k=1}^\infty 
    {T^{[h]}(x^k)\over k}\Big)
    \qquad\text{where}\quad
    T^{[h]}(x)=\sum_{n=1}^\infty T_n^{[h]}x^n
$$
\smallskip

{\bf 3.} To prove Proposition 7
we used
estimates
   $\log B_d < \log T_n^{[a\sqrt n]} = \log T_n - O(\sqrt n\,\log n)$
where $n=l(d)\sim \frac12 d^2$. 
A numerical experiment (using the above recurrent formula) 
allows us to surmise a better estimate 
$\log T_n^{[a\sqrt n]} = \log T_n - O(\log n)$.
\endremark

\head 2. Generalizations
\endhead

Some of the preceding results can be easily extended
to hypersurfaces in nonsingular varieties of any dimension.
To state these generalizations let us introduce
an additional notation.
For a very ample divisor $L$ on a nonsingular
variety $X$ of dimension $n$ denote by
$I(n,dL)$, respectively by $D(n,dL)$, the number of real nonsingular
hypersurfaces in the linear system $|dL|$ considered
up to isotopy, respectively up  to deformation equivalence.
As in the case of projective hypersurfaces, the number
$D(n,dL)$ is the number of connected components
of the complement of the corresponding discriminant
$\Delta(dL)\subset\vert dL\vert$
of the linear system $\vert dL\vert$. Let us call the
system {\it generic} if generic points of the discriminant
are the hypersurfaces $H\in\vert dL\vert$ with
a nondegenerate quadratic double point
which is the only its singular point.
Since for sufficiently large $d$ the projective image of $X$ by
$dL$ contains no linear subspaces, the system $\vert dL\vert$
is generic for large $d$ (see \cite{\refGH}).

\subhead      2.1. Hypersurfaces up to isotopy
\endsubhead
The lower
bound
in Proposition 1 extends to hypersurfaces
in any toric variety of any dimension.

\proclaim{ Proposition 10 } For any ample divisor $L$
on a toric variety $X$ of dimension $n$,
there is a constant $c=c(L,X)$ such that
$I(n,dL)\ge c\exp(d^n)$ for any $d$.
\endproclaim

%This extension of the lower
%bound from Proposition 1
%is proved in essentially the same manner
%as the one from the Proposition 1.
%     \mnote{Or: Vykinul objasnenie (i tak vse jasno)
%Kh: verojatno tyi prav}
%
% One notes that it is sufficient to treat
% the case when the Newton polyhedron is the cube,
% then, subdivides the cube with the edge
% $2m_k+6$ in $2^{n+1}$ cubes with the edge
% $m_k+2$ separated by slices of width.
%

The proof is essentially the same as that 
of the lower bound in Proposition 1.

%%%%%%%%%%%%%%%%%%%%%%%%%%%%%%%%%%%%%%%%%%%%%%%%%%%%%%%%%%%%%

\subhead     2.2. Hypersurfaces up to deformation
\endsubhead
The upper bound from Proposition 2 is
easily generalized to ample divisors on
nonsingular varieties of any dimension.
Via evident $I(n,dL)\le D(n,dL)$
it gives the best upper bound for $I(n,dL)$ we know.

\proclaim{Proposition 11} If $L$
is an ample divisor on a nonsingular variety of
dimension $n$, then, there 
is a constant $c'_n$ such that $$ \log I(n,dL)\le \log D(n,dL)\le
c'_nd^n\log d. 
%n.  
$$ \endproclaim

\demo{Proof}
As in the proof of Proposition 2, it is sufficient
to count the leading terms in the polynomial
expansion of the degree of the discriminant and
of the dimension of the linear system.
The first is counted through integrating the Euler
characteristic over a Lefschetz pencil: it gives
$a=(n+1)d^nL^n$. The second is counted by means of
the asymptotic Riemann-Roch formula: it gives
$b=\frac1{n!}d^nL^n$. Thus, it remains to
apply the Smith-Thom inequality, which gives
$a^b$ as the leading term of the bound.\qed
\enddemo

{\bf Remark.} There are two other approaches to bounding
from above $I$ and $D$. The first one
consists in constructing of
an algorithm recognizing the topological or $C^\infty$
type of a hypersurface and bounding the number of classes
by the complexity of the algorithm. Unfortunately,
we do not know algorithms which provide a better bound
than the one from Proposition 10.
The other one would work for hypersurfaces on toric
varieties and consists in studying of the
Gelfand-Kapranov-Zelvinsky secondary
polytopes, which are the Newton polytopes of the discriminants.
We collected some information concerning this approach in
the next subsection.

\subhead     2.3. Number of $T$-hypersurfaces
\endsubhead
Denote
the number of deformation
classes of $T$-hypersurfaces in $\vert dL\vert$
on a given
toric variety of dimension $n$
by $D'(n,dL)$.

\proclaim{Proposition 12}
There is a constant $c=c(n,L)$ such that
$\log D'(n,dL)\le cd^n.$
\endproclaim

The proof is similar to that of Proposition 3.

\subhead     2.4. Newton polytopes of the discriminants
\endsubhead
Let consider for simplicity the case of
hypersurfaces of projective space of dimension $n$.  In several
proofs above there appear the discriminant hypersurface
$\Delta(d)\subset P^N$ and its Newton polyhedron 
$Q(d)\subset\R^{n}$.
The
Betti numbers of $\Delta(d)$ bound the
Betti numbers of 
%$\RP^N\setminus\R\Delta(d)$.
$P^N_\R\setminus\Delta(d)_\R$.
The vertices of $Q(d)$
are in one-to-one correspondence with
convex integral triangulations of the simplex
$\{(i_1,\dots,i_n)\vert i_k\ge 0, \sum i_k\le d\}$
(see \cite{\refGKZ}).
In the proofs above, to bound the Betti numbers we used
only the degree of $\Delta(d)$ and to bound the number of
triangulations (and thus the number $D^T(n,d)$
of deformation classes of $T$-hypersurfaces)
the description, due to
Gelfand-Kapranov-Zelvinsky, of the coordinates of the vertices of
$Q(d)$.

Certainly, these two objects, $\Delta(d)$ and $Q(d)$,
contain more information which is still to be recovered.
For example, $\Delta(d)$ is a very singular variety,
which is not taken into account in our proofs. One could
expect to get a better bound on its total Betti number
through the number of critical values of the discriminant
polynomials, which is equivalent to bounding of the volume
of $Q(d)$. Unfortunately, as it is shown in \cite{\refOrevkovVol},
it is not the case. So, some finer analysis of
the topology of $\Delta(d)$ is needed.

It is worth noting that there is a related problem of
bounding the number of integral triangulations. 
Above we used already the observation that the number
of convex integral triangulations grows as
$\exp (cd^n)$. The nonconvex triangulations are related
somehow with the internal integral points of $Q(d)$
and one may expect that the number of nonconvex triangulations
has also the growth $\exp (cd^n)$ with, probably,
another constant $c$. For the moment such a result is proven
only for $n=2$, and
without appealing to $Q(d)$, see \cite{\refOrevkovJCTA}.

A natural question is {\it how many (asymptotically)
convex triangulations are among all integral triangulations?}
This is not clear even for $n=2$. As we already mentioned,
the both quantities are $\ewsim\exp(d^2)$ but do the coefficients
of $d^2$ coincide?
As a lower bound for the number $N_d$ of
all integral triangulations of a square
$d\times d$, one can use $(N_{k,m}^{1/km})^{d^2}$
where $N_{k,m}$ is the number of integral triangulations of
a rectangle $k\times m$ for some fixed $k$, $m$
(if $k$ and $m$ are fixed then $N_{k,m}$ can be computed explicitely).
In all the cases where we computed $N_{k,m}$, this bound was
much less than the upper bound $(4^4/3^3)^{d^2}$ for
the number convex triangulations
provided by the proof of Proposition 3.

Probably, the constant $4^4/3^3$ can be reduced as follows.
To get this constant, we estimated the number of vertices of
the secondary polytope
via the number of all integral points in the simplex containing it.
One can expect to obtain better estimates for the number
of vertices of a polytope if all the vertices are integral
(compare with 
\cite{\refArnold}, \cite{\refKS}, \cite{\refVershik}).

%%%%%%%%%%%%%%%%%%%%%%%%%%%%%%%%%%%%%%%%%%%%%%%%%%%%%%%%%%%%%
%%%%%%%%%%%%%%%%%%%%%%%%%%%%%%%%%%%%%%%%%%%%%%%%%%%%%%%%%%%%%%
%%%%%%%%%%%%%%%%%%%%%%%%%%%%%%%%%%%%%%%%%%%%%%%%%%%%%%%%%%%%%%
%%%%%%%%%%%%%%%%%%%%%%%%%%%%%%%%%%%%%%%%%%%%%%%%%%%%%%%%%%%%%%

\enddocument